\documentclass[11pt]{article}
\usepackage{amssymb}
\usepackage{}

\usepackage{authblk}
\usepackage{amsmath, amsthm, amssymb}
\usepackage{enumerate}
\usepackage{graphicx}
\usepackage{float}
\usepackage{hyperref}
\usepackage[margin=1in]{geometry}
\usepackage{comment}
\numberwithin{equation}{section}

\begin{document}

\theoremstyle{plain}
\newtheorem{theorem}{Theorem}[section]
\newtheorem{proposition}[theorem]{Proposition}
\newtheorem{lemma}[theorem]{Lemma}
\newtheorem{cor}[theorem]{Corollary}

\theoremstyle{definition}
\newtheorem{definition}[theorem]{Definition}
\newtheorem{remark}[theorem]{Remark}

\title{Shooting Method with Sign-Changing Nonlinearity}

\author{Ze Cheng$^{a}$ \footnote{ze.cheng@colorado.edu}
 and Congming Li$^{b,a}$  \footnote{congmingli@gmail.com}
 \footnote{Research partially supported by NSFC-11271166.} \\
$^a$ Department of Applied Mathematics, \\
University of Colorado Boulder, CO 80309, USA \\
$^b$ Department of Mathematics, INS and MOE-LSC, \\
Shanghai Jiao Tong University, Shanghai, China
}

\maketitle
\date{}

\begin{abstract}
In this paper, we study the existence of solution to a nonlinear system:
\begin{align}
  \left\{\begin{array}{cl}
      	-\Delta u_{i} = f_{i}(u) & \text{in } \mathbb{R}^n, \\
	       u_{i} > 0             & \text{in } \mathbb{R}^n, \, i = 1, 2,\cdots, L \\
        \end{array}
 \right.
\end{align}
for sign changing nonlinearities $f_i$'s. Recently, a degree theory approach to shooting method
for this broad class of problems is introduced in \cite{LiarXiv13} for nonnegative $f_i$'s.
However, many systems of nonlinear Sch\"odinger type involve interaction with undetermined sign.
Here, based on some new dynamic estimates, we  are able to extend the degree theory approach to
systems with sign-changing source terms.
\end{abstract}



\newpage

\section{Introduction}
Consider the systems
\begin{align}\label{original}
  \left\{\begin{array}{cl}
      	-\Delta u_{i} = f_{i}(u) & \text{in } \mathbb{R}^n, \\
	       u_{i} > 0             & \text{in } \mathbb{R}^n, \\
        \end{array}
 \right.
\end{align}
where $i = 1, 2,\cdots, L$. Denote $\mathbb{R}^{L}_{+} = \{u\in\mathbb{R}^L|u_i>0, \text{ for } i=1,\cdots,L\}$, and throughout this article $f = (f_1, f_2,\cdots, f_L)$
is assumed continuous in $\overline{\mathbb{R}^{L}_{+}}$
and locally Lipschitz continuous in $\mathbb{R}^{L}_{+}$.
The existence and nonexistence of positive solution of this type of systems has long been studied.
An important example (after reducing the order) is the Hardy-Littlewood-Sobolev type of system,
\begin{equation} \label{HLS}
\left\{ \begin{aligned}
         (-\triangle)^k u &= v^p & \text{in } \mathbb{R}^n,\\
         (-\triangle)^k v &= u^q & \text{in } \mathbb{R}^n,
                          \end{aligned} \right.
                          \end{equation}
of which the Lane-Emden system is a special case (k=1).
The Lane-Emden system possesses positive solutions in critical and supercritical cases,
i.e. $\frac{1}{p+1}+\frac{1}{q+1}\leq 1-\frac{2}{n}$
(cf. \cite{SZ-1998});
meanwhile, the system admits no radial positive solution in subcritical cases, i.e. $\frac{1}{p+1}+\frac{1}{q+1}> 1-\frac{2}{n}$,
by Mitidieri \cite{Mitidieri1996}.
The Lane-Emden conjecture says that the system admits no positive solution in subcritical cases.
In \cite{CL2009-DCDS}, W. Chen and C. Li show that under an integrability condition, \eqref{HLS} with $k=1$ is equivalent to a system of integral equations of which all positive solutions are radial.
This, together with Mitidieri's result, partially solves the Lane-Emden conjecture.
Also, Souplet solves the 4-dimension case of Lane-Emden conjecture in \cite{Souplet2009}. 
For more Liouville type of theorems, we refer readers to \cite{BrezisNirenberg1983, CGS, SG1981, NS} and the reference therein. 

Traditionally, one can use a concentrated-compactness argument (cf. \cite{Lions1985}) to show existence of solution to critical cases of Lane-Emden systems.
In \cite{LiarXiv13}, C. Li introduces a new approach to obtain existence of radial positive solution for system \eqref{original} in critical and supercritical cases, which connects shooting method with degree theory and surprisingly relates the question to the non-existence of solution to a corresponding Dirichlet problem.
For recent development of degree approach to shooting method, reader can check \cite{LiVillavert2013, Villavert2014}.

The conditions Li places on $f$ include that
$f$ needs to be positive (i.e. $f_i$ is positive for all $i=1,\cdots,L$)
 and satisfies a non-degenerate condition.
Although this covers a broad class of problems including critical and supercritical cases of \eqref{HLS}, there are cases of nonlinear Shr\"odinger type of systems that $f$ does not need to be always positive.
In this article, we are going to replace the condition $f$ being positive by $\sum f_i \geq 0$
where each $f_i$ can change sign.

Since we are looking for positive radial solution $u(x)=u(|x|)$ of \eqref{original}, the problem is equivalent to looking for global positive solution to the following ODE system,
\begin{align}\label{ode}
 \left\{ \begin{array}{cl}
u_{i}^{''}(r) + \dfrac{n-1}{r}u_{i}^{'}(r) = -f_{i}(u(r)) \\
u_{i}^{'}(0) = 0, \ u_{i}(0) = \alpha_i \text{ for } i = 1,2,\cdots, L.
\end{array}\right.
\end{align}
where $\alpha = (\alpha_1, \alpha_2, \cdots, \alpha_L)$
is positive (i.e. each $\alpha_i> 0$) initial value for $u$.

By classical ODE theory, this initial value problem \eqref{ode} has a unique solution $u(r,\alpha)$
for $r$ in some maximum interval.
Let $r_{\alpha} := \inf_{r \geq 0}\{r\in\mathbb{R} | u(r,\alpha)\in \partial\mathbb{R}^{L}_{+}\}$
(by $\partial\mathbb{R}^{L}_{+}$ we mean the boundary of $\mathbb{R}^{L}_{+}$, which sometimes we call ``the wall"),
so $r=r_{\alpha}$ is where $u(r,\alpha)$ touches the wall for the first time.
There are two cases, case 1: $r_{\alpha}=\infty$, $u$ never hits wall,
then $u(r,\alpha)$ is a radial positive solution to \eqref{original}, and we are done; case 2: $r_{\alpha}<\infty$, $u$ hits wall in finite time, and
we will show that the existence of solution of \eqref{original} surprisingly depends on the non-existence of solution to its corresponding Dirichlet problem \eqref{Dirichlet} in the below,

\begin{align}\label{Dirichlet}
  \left\{\begin{array}{cl}
      	-\Delta u_{i} = f_{i}(u) & \text{in } B_R, \\
	       u_{i} > 0             & \text{in } B_R, \\
           u_i=0                 & \text{on } \partial B_R,
        \end{array}
 \right.
\end{align}
where $B_R=B_R(0)$ for any $R>0$ and $i = 1,2,\cdots, L$. 

It is well known that the nonexistence of solution 
to this critical and supercritical Dirichlet problems \eqref{Dirichlet} 
can be derived by Rellich-Poho\v{z}aev type of identities.
In fact, it belongs to a general variational problem, where
a ball of any radius can be replaced by a bounded star-shaped domain.
For this problem, we refer the reader to Pucci and Serrin's paper \cite{PS-Ind}.
In the last part of this paper, we will show the local nonexistence to prove global existence
to some example systems.

So, here is our main result,
\begin{theorem}\label{existencetheorem}
Given the nonexistence of solution to system \eqref{Dirichlet} for all $R>0$, the system \eqref{original} admits a radially symmetric solution of class $C^{2,\alpha}(\mathbb{R}^{n})$ with $0<\alpha<1$, if $f= (f_1(u), f_2(u),\cdots, f_L(u)):\mathbb{R}^L\rightarrow \mathbb{R}^L$ satisfies the following assumptions:
\begin{enumerate}
  \item $f$ is continuous in $\overline{\mathbb{R}^{L}_{+}}$
and locally Lipschitz continuous in $\mathbb{R}^{L}_{+}$, and furthermore,
  \begin{align}\label{decayAssumption}
        \sum_{i=1}^L f_i(u) \geq 0  \text{ in } \mathbb{R}_+^L;
  \end{align}
  \item If $\alpha \in \partial\mathbb{R}^{L}_{+}$ and $\alpha\neq 0$, i.e., for some permutation $(i_1,\cdots,i_L)$, $\alpha_{i_1}=\cdots=\alpha_{i_m}=0$, $\alpha_{i_{m+1}},\cdots, \alpha_{i_L}>0$ where $m$ is an integer in $(0,L)$, then $\exists \delta_0 = \delta_0(\alpha) > 0$ such that for $\beta \in \mathbb{R}^{L}_{+}$ and $|\beta-\alpha|<\delta_0$,
  \begin{align}\label{ControlInequality}
    \sum_{j=m+1}^L |f_{i_j}(\beta)| \leq C(\alpha) \sum_{j=1}^m f_{i_j}(\beta),
  \end{align}
  where $C$ is a non-negative constant that depends only on $\alpha$.
\end{enumerate}
\end{theorem}

\begin{remark}\label{assumptionRemark}
Notice that under assumptions above, sign-changing source term $f$ is allowed. It is known that if $f_i$'s stay positive, there are many nice properties that we can use,  c.f. \cite{LiarXiv13, LiVillavert2013, Villavert2014}. If $f_i$ changes sign, those properties are lost, which leads to estimates failing. Our work here is to derive dynamic estimate \eqref{dynamicEstimate1} and \eqref{dynamicEstimate2} under assumptions above, such that the degree theory approach to shooting method is applicable to show existence of solution with sign-changing $f$.

\end{remark}

The paper is organized as follows. 
In section 2, we prove the main result.
In section 3, we will show nonexistence of solution to the corresponding Dirichlet problems
of some example systems.

\section{Proof of main result}
In this section, we will first define a target map and prove its continuity. Then we apply degree theory to prove theorem \ref{existencetheorem}.

\subsection{Target map}
For any real number $a>0$, let $\Sigma_a=\{\alpha\in\overline{\mathbb{R}^{L}_{+}}|\sum_{i=1}^L \alpha_i=a\}$, and $B_a= \{ \alpha \in \partial \mathbb{R}_+^L  {\bf |}
\sum_{i=1}^L \alpha_i \leq a\}$. Recall that for positive $\alpha$ (i.e. every $\alpha_i>0$) we define $r_{\alpha} = \inf_{r \geq 0}\{r\in\mathbb{R} | u(r,\alpha)\in \partial\mathbb{R}^{L}_{+}\}$. As mentioned before, we can assume $r_{\alpha}<\infty$ since if $r_{\alpha}=\infty$ we get a solution to \eqref{original}. Then we define a target map on a initial data as following,
\begin{definition}\label{psi}
Let $u(r,\alpha)$ be a solution to \eqref{ode} with initial value $\alpha\in\Sigma_a$, we define a map
$\psi : \Sigma_a \rightarrow \partial\mathbb{R}^{L}_{+}$, such that
\begin{align}
    \psi(\alpha) = \left\{\begin{array}{ll}
      	 u(r_{\alpha},\alpha) & \alpha\in\mathbb{R}^{L}_{+}, \\
	       \alpha            & \alpha\in\partial\mathbb{R}^{L}_{+}.
        \end{array}
 \right.
\end{align}

\end{definition}

Here we sketch Li's degree theory approach for shooting method 
as follows. 
Fix any real number $a>0$,
and assume that for any initial value $\alpha\in\Sigma_a$
no global positive solution to \eqref{ode} exists
(i.e. $r_{\alpha}<\infty$), so we can define a target map. Hence,
{\bf step 1},  we show that, under some suitable assumptions on $f$,
the target map $\psi$ is continuous from $\Sigma_a$ to $B_a$;
{\bf step 2}, by degree theory we show that $\psi$ is onto,
therefore $\exists \alpha_0\in\Sigma_a$ such that $\psi(\alpha_0)=u(r_{\alpha_0},\alpha_0)=0$;
{\bf step 3}, note that by assumption $r_{\alpha_0}<\infty$, $u(r,\alpha_0)$
for $r\in[0,r_{\alpha_0}]$ is a solution to the Dirichlet problem \eqref{Dirichlet}, 
which makes a contradiction 
if we assume that system \eqref{Dirichlet} admits no solution for any $R>0$.

In what follows, we assume \eqref{ode} admits no global positive solution, i.e. $r_{\alpha}<\infty$. We first show that under our assumptions \eqref{decayAssumption} and \eqref{ControlInequality} on $f$, the behavior of $u$ turns regular, such that $\psi$ is continuous. 
\begin{lemma}\label{continuityLemma}
For any real number $a>0$, let
\begin{align*}
    \Sigma_a=\{\alpha\in\overline{\mathbb{R}^{L}_{+}}|\sum_{i=1}^L \alpha_i=a\} \text{ and } B_a= \{ \alpha \in \partial \mathbb{R}_+^L | \sum_{i=1}^L \alpha_i \leq a\}.
\end{align*}
The target map $\psi$ defined in definition \ref{psi} is a continuous map from $\Sigma_a$ to $B_a$ if $f$ satisfies assumptions \eqref{decayAssumption} and \eqref{ControlInequality}.
\end{lemma}



Proof. 
To see that $\psi$ maps $\Sigma_a$ to $B_a$, we need to notice that by assumption \eqref{decayAssumption} $\Sigma_{i=1}^L f_i\geq 0$, so we solve from the ODE system \eqref{ode} and get
\begin{align*}
    \Sigma_{i=1}^L u_i(r,\alpha) &= \Sigma_{i=1}^L \alpha_i - \Sigma_{i=1}^L\int_0^r\int_0^s (\frac{\tau}{s})^{n-1} f_i(u(\tau)) d\tau ds \\
    &\leq \Sigma_{i=1}^L \alpha_i,
\end{align*}
for $r\in[0,r_\alpha]$. Therefore, $\psi(\alpha)\in B_a$.

Next, we will show that $\psi$ is continuous on $\Sigma_a$. Fix any $\overline{\alpha}\in\Sigma_a$, then $\overline{\alpha}$ lies on the boundary of $\mathbb{R}^{L}_+$ or in $\mathbb{R}^{L}_+$ ($\overline{\alpha}\neq 0$ since $a>0$), and we will prove for these two cases.

Case 1. $\overline{\alpha} \in \partial\mathbb{R}^{L}_+$. 

Suppose $\overline{\alpha}_{i_1}=\cdots=\overline{\alpha}_{i_m}=0$, and $\overline{\alpha}_{i_{m+1}},\cdots, \overline{\alpha}_{i_L}>0$, for some integer $m$ that $0<m<L$.

By the second assumption \eqref{ControlInequality}, $\exists \delta_0>0$, such that for $\alpha \in \mathbb{R}^{L}_{+}$ satisfying $|\alpha-\overline{\alpha}|<\delta_0$ we have
\begin{align}\label{environmentControl}
    \sum_{j=m+1}^L |f_{i_j}(\alpha)| \leq C(\overline{\alpha}) \sum_{1\leq j\leq m} f_{i_j}(\alpha).
\end{align}
Notice that we can choose $C=C(\overline{\alpha})\geq 1$ in \eqref{environmentControl} if $\sum_{1\leq j\leq m} f_{i_j}(\alpha)\geq 0$. Indeed,
if $C=0$, the term on the left of the inequality \eqref{environmentControl} is zero, and due to the first assumption \eqref{decayAssumption}, $\sum_{1\leq j\leq m} f_{i_j}(\alpha)\geq 0$. If $C>0$, then $\sum_{1\leq j\leq m} f_{i_j}(\alpha)\geq 0$ obviously. So, we choose $C\geq 1$.

For this $\delta_0$ and $C$, we {\bf{claim}} that:

\emph{For $\alpha\in\Sigma_a$, fix any $\delta < \frac{\delta_0}{2(3+L)C}$, and if $|\alpha-\overline{\alpha}|\leq \delta$, then for $r\in[0,r_{\alpha}]$
\begin{align}\label{dynamicEstimate1}
    |u(r,\alpha)-\overline{\alpha}| \leq 2(3+L)C\delta < \delta_0.
\end{align}
}

As we will see in the following proof, \eqref{dynamicEstimate1} is a dynamic estimate, 
in the sense that if $|u(r,\alpha)-\overline{\alpha}|<\delta_0$ 
with $r\in[0, a_1)\subset[0,r_{\alpha}]$ for some $a_1>0$, 
then by \eqref{ControlInequality} we have
\begin{align}
		\sum_{j=m+1}^L |f_{i_j}(u(r,\alpha)| \leq C(\overline{\alpha}) \sum_{1\leq j\leq m} f_{i_j}(u(r,\alpha)),
\end{align}
and this control enable us to push the range of $r$ in 
$|u(r,\alpha)-\overline{\alpha}|<\delta_0$
further than $a_1$ and up to $r_{\alpha}$. 

Suppose the claim not true, then there exists $\alpha_0 \in\mathbb{R}^{L}_{+}$
satisfying $|\alpha_0-\overline{\alpha}|\leq\delta$ and $a_1\in(0,r_{\alpha_0})$ such that the equality of \eqref{dynamicEstimate1} holds at $r=a_1$ for the first time, i.e.
\begin{align}\label{firstTimeEquality}
    |u(r,\alpha_0)-\overline{\alpha}|\left\{\begin{array}{ll}
      	 < 2(3+L)C\delta, & \text{if } r<a_1, \\
	     = 2(3+L)C\delta, & \text{if } r=a_1.
        \end{array}
 \right.
\end{align}
For $r\in(0,a_1)$ we have,
\begin{align}
    |u(r,\alpha_0)-\overline{\alpha}|
    &\leq |u(r,\alpha_0)-\alpha_0| + |\alpha_0-\overline{\alpha}|\\
    &\leq \sum_{j=1}^m|u_{i_j}(r,\alpha_0)-\alpha_{0i_j}|+\sum_{j=m+1}^L|u_{i_j}(r,\alpha_0)-\alpha_{0i_j}| + |\alpha_0-\overline{\alpha}|. \label{threeTerms}
\end{align}
So, to estimate the second term of \eqref{threeTerms}, we solve from \eqref{ode} and get
\begin{align*}
    \sum_{j=m+1}^L|u_{i_j}(r,\alpha_0)-\alpha_{0i_j}|
        &= \sum_{j=m+1}^L|\int_0^r\int_0^s (\frac{\tau}{s})^{n-1} f_{i_j}(u(\tau)) d\tau ds| \\
\end{align*}        
Notice that since $C\geq 1$, $2(3+L)C\delta<\delta_0$, therefore assumption \eqref{ControlInequality} can be applied on $u(r,\alpha_0)$ for $r\in(0,a_1)$,
\begin{align*}
    \sum_{j=m+1}^L|u_{i_j}(r,\alpha_0)-\alpha_{0i_j}| 
        &\leq \sum_{j=m+1}^L \int_0^r\int_0^s |(\frac{\tau}{s})^{n-1} f_{i_j}(u(\tau))| d\tau ds \\
        &\leq C \int_0^r\int_0^s (\frac{\tau}{s})^{n-1} \sum_{j=1}^m f_{i_j}(u(\tau)) d\tau ds \\
        &= C \sum_{j=1}^m (\alpha_{0i_j}-u_{i_j}(r,\alpha_0)),
\end{align*}

The first term of \eqref{threeTerms} can be estimated by
\begin{align*}
    \sum_{j=1}^m|u_{i_j}(r,\alpha_0)-\alpha_{0i_j}|
        &= \sum_{j=1}^m(\alpha_{0i_j}-u_{i_j}(r,\alpha_0))^+ + \sum_{j=1}^m(\alpha_{0i_j}-u_{i_j}(r,\alpha_0))^- \\
        &\leq 2\sum_{j=1}^m(\alpha_{0i_j}-u_{i_j}(r,\alpha_0))^+.
\end{align*}
To see the inequality above, one needs to notice that due to \eqref{ControlInequality}, $\sum_{j=1}^m f_{i_j}(u)\geq 0$, so $\sum_{j=1}^m u_{i_j}(r,\alpha_0)$ is monotone decreasing on $[0,r_{\alpha_0}]$. So,
\begin{align*}
    \sum_{j=1}^m(\alpha_{0i_j} - u_{i_j}(r,\alpha_0))
            = \sum_{j=1}^m(\alpha_{0i_j}-u_{i_j}(r,\alpha_0))^+ - \sum_{j=1}^m(\alpha_{0i_j}-u_{i_j}(r,\alpha_0))^- \geq 0.
\end{align*}

The last term of \eqref{threeTerms} is bounded by $\delta$, and notice that $u_i(r,\alpha_0)>0$, $i=1,\cdots, L$ for $r\in(0,r_{\alpha_0})$, so we get
\begin{align*}
    |u(r,\alpha_0)-\overline{\alpha}|
    &\leq 2\sum_{j=1}^m(\alpha_{0i_j}-u_{i_j}(r,\alpha_0))^+ + C \sum_{j=1}^m (\alpha_{0i_j}-u_{i_j}(r,\alpha_0)) + \delta\\
    &\leq (2+C) \sum_{j=1}^m \alpha_{0i_j} + \delta \\
    &\leq (2+C)L|\alpha_0-\overline{\alpha}| + \delta \\
    &\leq (3+C)L\delta,
\end{align*}
where $C$ is the same $C$ in \eqref{ControlInequality}. Then we get a contradiction with \eqref{firstTimeEquality} by taking $r=a_1$ in the above. Hence the claim is proved.

Notice that the claim and estimate \eqref{dynamicEstimate1} implies the continuity of $\psi$ at $\overline{\alpha}$, therefore, we have proved case 1.

Case 2. $\overline{\alpha} \in \mathbb{R}^{L}_+$.

As above we assume $r_{\overline{\alpha}}<\infty$, and obviously  $r_{\overline{\alpha}}>0$. Let's assume at $r=r_{\overline{\alpha}}$, for some integer $m$ that $0< m\leq L$ ($m>0$ because $\psi(\overline{\alpha})=u(r_{\overline{\alpha}},\overline{\alpha})\in\partial\mathbb{R}^L_+$, i.e. $u$ touches wall at $r=r_{\overline{\alpha}}$), so suppose $u_{i_1}(r_{\overline{\alpha}},\overline{\alpha})=\cdots=u_{i_m}(r_{\overline{\alpha}},\overline{\alpha})=0$, and $u_{i_{m+1}}(r_{\overline{\alpha}},\overline{\alpha}),\cdots,u_{i_L}(r_{\overline{\alpha}},\overline{\alpha})>0$. Let $\overline{\omega}(r)=\sum_{j=1}^m u_{i_j}(r,\overline{\alpha})$, and we claim that $\overline{\omega}'(r_{\overline{\alpha}})<0$.

By the continuity of $u(r,\overline{\alpha})$ with respect to $r$, $\exists \delta_1>0$ such that for $r\in(r_{\overline{\alpha}}-\delta_1,r_{\overline{\alpha}}]$,
\begin{align}
    |u(r,\overline{\alpha})-\psi(\overline{\alpha})|<\delta_0.
\end{align}
If $m<L$, the assumption \eqref{ControlInequality} is taking effect, and therefore
 $\sum_{j=1}^m f_{i_j} \geq 0$ for $r\in(r_{\overline{\alpha}}-\delta_1,r_{\overline{\alpha}}]$;
if $m=L$ then by the assumption \eqref{decayAssumption},
 we also have $\sum_{j=1}^m f_{i_j} \geq 0$ for $r\in(r_{\overline{\alpha}}-\delta_1,r_{\overline{\alpha}}]$.
So, in $(r_{\overline{\alpha}}-\delta_1,r_{\overline{\alpha}}]$
\begin{align}\label{dynamicEstimate2}
    -\frac{1}{r^{n-1}}(r^{n-1}\overline{\omega}'(r))'=\sum_{j=1}^m f_{i_j} \geq 0.
\end{align}
Also, since $\overline{\omega}(r)>0$ for $r<r_{\overline{\alpha}}$ and $\overline{\omega}(r_{\overline{\alpha}})=0$,
there must exist $r_0\in(r_{\overline{\alpha}}-\delta_1,r_{\overline{\alpha}})$, such that $\overline{\omega}'(r_0)<0$.
So, for $r\in[r_0,r_{\overline{\alpha}}]$,
\begin{align}
    \overline{\omega}'(r) = (\frac{r_0}{r})^{n-1}\overline{\omega}'(r_0)
                    -\int_{r_0}^r(\frac{\tau}{r})^{n-1}\sum_{j=1}^m f_{i_j}(u(\tau))d\tau<0.
\end{align}
This proves the claim $\overline{\omega}'(r_{\overline{\alpha}})<0$.
Then there exists $l_0\in\{1,\cdots, m\}$
such that $u_{i_{l_0}}'(r_{\overline{\alpha}},\overline{\alpha})<0$.
Therefore, combining with the fact $u_{i_{l_0}}(r_{\overline{\alpha}},\overline{\alpha})=0$, we see $u_{i_{l_0}}$ crosses the wall with a non-zero slope, i.e. there is a transversality at $r=r_{\overline{\alpha}}$, and by classical ODE stability theory (the continuous dependence on initial value) $\psi$ is continuous at $\overline{\alpha}$. $\Box$

\subsection{Application of degree theory}
Now, let's recall some results in degree theory (in particular, the treatment modified by P.Lax, c.f. \cite{Nirenberg2001}). Consider $C^{\infty}$ oriented manifolds $X_0,Y$ of dimension $n$ (all manifolds are assumed to be paracompact) and an open subset $X\subseteq X_0$ with compact closure. For convenience, write $dy^1\wedge\cdots\wedge dy^n = dy$. Then for a $C^1$ map $\phi:X\rightarrow Y$, the degree is defined as following:
\begin{definition}\label{degree}
Let $\mu=f(y)dy$ be a $C^{\infty}$ $n$-form with support contained in a coordinate patch $\Omega$ of $y_0$ and lying in $Y\setminus\{\phi(\partial X)\}$ such that $\int_Y \mu=1$; set
\begin{align}
    deg(\phi,X,y_0)=\int_X \mu\circ\phi.
\end{align}
\end{definition}

Here are some properties of degree which we will refer to in our proof,
\begin{proposition}
For $y_1$ close to $y_0$, $deg(\phi,X,y_0)=deg(\phi,X,y_1)$.
\end{proposition}
It follows that the degree of a mapping is constant on any connected component $C$ of $Y\setminus \{\phi(\partial X)\}$, and we can write degree as $deg(\phi,X,C)$.

\begin{proposition}\label{ontoProperty}
If $y_0\notin\phi(\overline{X})$, then $deg(\phi,X,y_0)=0$.
\end{proposition}
As a result, if $deg(\phi,X,y_0)\neq0$, then $y_0\in\phi(\overline{X})$.

An important property of degree is that, the notion can be extended to maps $\phi$ which are merely continuous, since we can approximate such $\phi$ by $C^1$ maps $\phi_n$ (see Property 1.5.3 in \cite{Nirenberg2001}). Also, degree is homotopy invariant, which enables us to define degree for continuous map $\eta:\partial X\rightarrow \mathbb{R}^n\setminus y_0$ (see Property 1.5.4 in \cite{Nirenberg2001}), which leads to the following theorem (see Property 1.5.5 in \cite{Nirenberg2001}),

\begin{theorem}\label{homotopyInvariance}
$deg(\eta, X, y_0)$ depends only on the homotopy class of $\eta: \partial X\rightarrow \mathbb{R}^n\setminus y_0$.
\end{theorem}
Now we are prepared to prove the existence of solution to \eqref{original}.\\

\noindent{\bf Proof of theorem \ref{existencetheorem}. }
For any fixed real number $a>0$, assume that \eqref{ode} admits no global positive solution with any initial value $\alpha\in\Sigma_a$, so $r_{\alpha}<\infty$, and then we can define a target map $\psi$ by \eqref{psi}.

Recall that $\Sigma_a=\{\alpha\in\overline{\mathbb{R}^{L}_{+}}|\sum_{i=1}^L \alpha_i=a\}$ and $B_a= \{ \alpha \in \partial \mathbb{R}_+^L  {\bf |}
\sum_{i=1,\cdots, L} \alpha_i \leq a\}$, and by lemma \ref{continuityLemma} $\psi$ is a continuous maps from $ \Sigma_a \longrightarrow B_a$.

Let $\pi(\alpha)=\alpha+\frac{1}{L}(a- \displaystyle\sum_{i=1,\cdots, L} \alpha_i)(1,\cdots,1): B_a \longrightarrow \Sigma_a$,
then $\pi$ is continuous with a continuous inverse $\pi^{-1}(\alpha)=\alpha-(\displaystyle\min_{i=1,\cdots, L} \alpha_i)(1,\cdots,1): \Sigma_a \longrightarrow B_a$.

The map: $\phi=\pi \circ \psi: \Sigma_a \longrightarrow \Sigma_a$ is continuous and $\phi(\alpha)=\alpha$
on $\partial\Sigma_a$. Let $\eta=id$ (the identity map) and $X=\Sigma_a\setminus \partial\Sigma_a$,
and then by theorem \ref{homotopyInvariance} we have $deg(\phi, X, \alpha)=deg(\eta, X, \alpha)=1$ for any
$\alpha \in \Sigma_a\setminus\partial\Sigma_a$. By property \ref{ontoProperty}, $\phi$ is onto, which implies that $\psi$ is also onto. this shows that there exists an $\alpha_0 \in \Sigma_a$ such that $\psi(\alpha_0)=0$.

Since we assume that system \eqref{Dirichlet} admits no solution, $r_{\alpha_0}$ corresponding to this $\alpha_0$ cannot be finite, a contradiction. This completes the proof of theorem \ref{existencetheorem}. $\Box$

\section{Examples}
One of the simplest systems is that $f\equiv 0$ in \eqref{original}, then $u\equiv C$ for some constant vector $C$ is a solution. Let us point out that since $f$ is not positive, this trivial system is not included in previous results (cf. \cite{LiarXiv13, LiVillavert2013, Villavert2014}), but it is included in our cases. In this section, we will show the existence of solution to some non-trivial systems. In the view of theorem \ref{existencetheorem}, we only need to show their corresponding Dirichlet problems admit no solution, and verify its source term satisfy our assumptions in theorem \ref{existencetheorem}.

\subsection{Sign-changing source terms}
Here we give a simple but non-trivial example of sign-changing source terms system, and we believe there are many other non-linear Schr\"odinger type of systems with sign-changing source terms which our method can be applied to.
Consider the following system,
\begin{align}\label{signChanging}
  \left\{\begin{array}{ll}
    -\Delta u = v^p-u^p,  \\
    -\Delta v = u^p,    \\
	\ \ u,v > 0, \\
        \end{array}
\right.\text{in } \mathbb{R}^n,
\end{align}
and its corresponding Dirichlet problem,
\begin{align}\label{signChangingDirichlet}
  \left\{\begin{array}{ll}
    -\Delta u = v^p-u^p & \text{in } B,  \\
    -\Delta v = u^p   & \text{in } B,  \\
	\ \ u,v > 0 & \text{in } B, \\
    \ \ u,v=0 & \text{on } \partial B,\\
        \end{array}
\right.
\end{align}
where $B=B_R(0)\subset \mathbb{R}^n$ for any $R>0$. We have
\begin{theorem}\label{signChangingExistenceThm}
If $p\geq\frac{n+2}{n-2}$, then \eqref{signChanging} admits radial positive solution.
\end{theorem}

Again the proof relies on the non-existence of solution to \eqref{signChangingDirichlet}. The non-existence is obtained by computing Rellich-Poho\v{z}aev type identity.

\begin{lemma}\label{signChangingDirichletNonEx}
If $p\geq\frac{n+2}{n-2}$, then system \eqref{signChangingDirichlet} admits no solution for any $R>0$.
\end{lemma}
Proof. Suppose there exists a positive solution $(u,v)$.

{\bf{Step 1.}} Claim the following identity,
\begin{align}\label{pohozaevCross}
\int_B \Delta u(x\cdot \nabla v) + \Delta v(x\cdot\nabla u) - (n-2)\nabla u\cdot \nabla v dx= \int_{\partial B} (x\cdot \nu)(\frac{\partial u}{\partial \nu}\frac{\partial v}{\partial \nu}) d\sigma>0,
\end{align}
where $\nu$ is the outward normal. The calculation of the above identity  usually goes as
\begin{align*}
	\int_B \Delta u(x\cdot \nabla v) dx &= \int_B \text{div} (\nabla u (x\cdot \nabla v)) - \nabla u\cdot \nabla(x\cdot \nabla v) dx \\
 				&= \int_{\partial B} (\nabla u \cdot \nu)(x\cdot \nabla v) d\sigma- \int_B ( \nabla u \cdot \nabla v + x_j \partial_i u \partial^2_{ij} v) dx.
\end{align*}
Then do the same to the second integrand in \eqref{pohozaevCross} and sum it with the above one to get
\begin{align*}
	&\int_B \Delta u(x\cdot \nabla v) + \Delta v(x\cdot\nabla u) dx \\
	&= \int_{\partial B} \frac{\partial u}{\partial \nu}(x\cdot \nabla v) +\frac{\partial v}{\partial \nu}(x\cdot \nabla u) d\sigma - \int_B (2 \nabla u \cdot \nabla v + x_j \partial_i u \partial^2_{ij} v+x_j \partial_i v \partial^2_{ij} u) dx \\
	&= \int_{\partial B} \frac{\partial u}{\partial \nu}(x\cdot \nabla v) +\frac{\partial v}{\partial \nu}(x\cdot \nabla u) d\sigma  - \int_B (2 \nabla u \cdot \nabla v + x\cdot\nabla(\nabla u\cdot\nabla v)) dx \\
	&= \int_{\partial B} \frac{\partial u}{\partial \nu}(x\cdot \nabla v) +\frac{\partial v}{\partial \nu}(x\cdot \nabla u) d\sigma  - \int_B 2 \nabla u \cdot \nabla v dx - \int_{\partial B} x\cdot \nu (\nabla u\cdot\nabla v)) d\sigma + n\int_B \nabla u\cdot\nabla v dx.
\end{align*}
Notice the fact that $x=|x|\nu$ on $\partial B$, and $\nabla u = \frac{\partial u}{\partial \nu} \nu$ and $\nabla v = \frac{\partial v}{\partial \nu} \nu$ due to $u,v=0$ on $\partial B$, and after rearrangement we have the identity \eqref{pohozaevCross}. Also, by Hopf's Lemma $\frac{\partial u}{\partial \nu}<0$ and $\frac{\partial v}{\partial \nu}<0$, so we have RHS of \eqref{pohozaevCross} $>0$.

Now, since
\begin{align*}
	-\int_B v\Delta u dx = \int_B \nabla u\cdot \nabla v dx = -\int_B u\Delta v dx,
\end{align*}
for $\theta\in[0,1]$ we have
\begin{align}\label{energyCross}
	\int_B \nabla u\cdot \nabla v dx = - \int_B \theta v\Delta u + (1-\theta) u\Delta vdx.
\end{align}
Also, notice that in \eqref{pohozaevCross} if we let $u=v$, then
\begin{align}\label{pohozaevSingle}
	\int_B \Delta u (x\cdot \nabla u)  dx= \int_B  \frac{n-2}{2}|\nabla u|^2 dx + \frac 1 2 \int_{\partial B} x\cdot \nu |\nabla u|^2 d\sigma>0. 
\end{align}

{\bf{Step 2.}} Merge the source terms into \eqref{pohozaevCross} (we replace $-u^p$ in the first source term by $\Delta v$, i.e. $-\Delta u=v^p+\Delta v$ and $-\Delta v = u^p$), and by \eqref{energyCross} we have
\begin{align*}
	&  \text{LHS of \eqref{pohozaevCross}} = \int_B -(v^p+\Delta v)(x\cdot \nabla v) - u^p(x\cdot \nabla u) - (n-2)\nabla u\cdot \nabla v dx \\
	&=\int_B -x\cdot \nabla\left( \frac{v^{p+1}}{p+1}+\frac{u^{p+1}}{p+1}\right)  - \Delta v (x\cdot \nabla v) 
	+ (n-2)\left(  \theta v\Delta u + (1-\theta) u\Delta v\right) dx \\
	&=\int_B -x\cdot \nabla\left( \frac{v^{p+1}}{p+1}+\frac{u^{p+1}}{p+1}\right)  - \Delta v (x\cdot \nabla v) 
		+ (n-2)\left(  -\theta (v^{p+1}+v\Delta v) - (1-\theta) u^{p+1}\right) dx \\
	&= \int_B \left\lbrace \left( \frac {n} {p+1} - (n-2)\theta\right)  v^{p+1} + \left( \frac {n} {p+1} - (n-2)(1-\theta)\right)  u^{p+1}  \right\rbrace dx \\
	&\quad - \int_B \Delta v (x\cdot \nabla v) + (n-2)\theta v\Delta v dx	 
\end{align*}
So, by \eqref{pohozaevSingle}, \eqref{pohozaevCross} becomes 
\begin{align*}
	&\int_B \left\lbrace \left( \frac {n} {p+1} - (n-2)\theta\right)  v^{p+1} + \left( \frac {n} {p+1} - (n-2)(1-\theta)\right)  u^{p+1}  \right\rbrace dx \\
	&= \int_B \Delta v (x\cdot \nabla v) + (n-2)\theta v\Delta v dx	 +\int_{\partial B} (x\cdot \nu)(\frac{\partial u}{\partial \nu}\frac{\partial v}{\partial \nu}) d\sigma \\
	&= \int_B \frac{n-2}{2}|\nabla v|^2 dx+\frac 1 2 \int_{\partial B} x\cdot \nu |\nabla u|^2 d\sigma 
	-\int_B\theta(n-2)|\nabla v|^2 dx
	+ \int_{\partial B} (x\cdot \nu)(\frac{\partial u}{\partial \nu}\frac{\partial v}{\partial \nu}) d\sigma \\
	&= \int_B (1-2\theta)\frac{n-2}{2}|\nabla v|^2 dx+\frac 1 2 \int_{\partial B} x\cdot \nu |\nabla u|^2 d\sigma 
		+ \int_{\partial B} (x\cdot \nu)(\frac{\partial u}{\partial \nu}\frac{\partial v}{\partial \nu}) d\sigma.
\end{align*}
Take $\theta=\frac 1 2$, and we have
\begin{align}\label{pohozaevEx1}
\int_B  \left( \frac {n} {p+1} - \frac{n-2}{2}\right)  \left( v^{p+1} + u^{p+1}\right) dx
= \frac 1 2 \int_{\partial B} x\cdot \nu |\nabla u|^2 d\sigma 
		+ \int_{\partial B} (x\cdot \nu)(\frac{\partial u}{\partial \nu}\frac{\partial v}{\partial \nu}) d\sigma >0.
\end{align}
So, by assumption $p\geq \frac{n+2}{n-2}$ the LHS of the above identity $\leq 0$, a contradiction.

$\Box$ \\

\noindent {\bf{Proof of Theorem \ref{signChangingExistenceThm}.}}
Directly we see that \eqref{signChanging} satisfies our first main assumption \eqref{decayAssumption}.
To see \eqref{ControlInequality} is satisfied, just notice that if $u=0$ then $|f_2|=0\leq f_1=v^p$, and if $v=0$ then $|f_1|=u^p\leq f_2$. Then for a neighborhood of such $(u,v)$ (i.e. $uv=0$), \eqref{ControlInequality} holds.

So, combined with Lemma \ref{signChangingDirichletNonEx}, Theorem \ref{signChangingExistenceThm} follows from Theorem \ref{existencetheorem}. \\
$\Box$

\begin{remark}
Actually system \eqref{signChanging} can be solved as follows. Suppose $-\Delta w = w^p$, let $u=\lambda w$ and $v=\nu w$, then we can find suitable $\lambda,\nu$ such that $(u,v)$ is a solution. 

So, to manifest a nontrivial application of Theorem \ref{existencetheorem} we consider a very similar system, 
\begin{align}\label{example2}
  \left\{\begin{array}{ll}
    -\Delta u = v^p +v^q -u^p, \\
    -\Delta v = u^p, \\
	\ \ u,v > 0,   \\
        \end{array}
\right. \text{in } \mathbb{R}^n,
\end{align}
where $p\neq q$ and $p,q\geq \frac{n+2}{n-2}$.
To show this system admits a solution, we only need to show the nonexistence of solution to corresponding Dirichlet problem (source terms satisfy assumptions \eqref{decayAssumption} and \eqref{ControlInequality}, see similar proof of Theorem \ref{signChangingExistenceThm}). Then the first step is exactly the same as the proof of Lemma \ref{signChangingDirichletNonEx}. In the second step, we use $-\Delta u = v^p+v^q+\Delta v$ and $-\Delta v=u^p$ to merge the source term, and taking $\theta=\frac 1 2$ \eqref{pohozaevEx1} becomes 
\begin{align}
\begin{array}{cc}
\int_B  \left( \frac {n} {p+1} - \frac{n-2}{2}\right)  \left( v^{p+1} + u^{p+1}\right) +\left( \frac{n}{q+1} -\frac{n-2}{2}\right)v^{q+1} dx  \quad\\
= \quad \frac 1 2 \int_{\partial B} x\cdot \nu |\nabla u|^2 d\sigma 
		+ \int_{\partial B} (x\cdot \nu)(\frac{\partial u}{\partial \nu}\frac{\partial v}{\partial \nu}) d\sigma >0.
\end{array}
\end{align}
Then the nonexistence to Dirichlet problem follows.
\end{remark}

\subsection{Conservative source terms}
In this section, we consider a system that $f$ has a potential function $F$, i.e. $f=\nabla F$, and $F(0)=0$.

{\bf{Type I.}} Consider the following system,s
\begin{align}\label{systemPotential}
  \left\{\begin{array}{cl}
      	-\Delta u_{i} = \frac{\partial F}{\partial u_i},  \\
	       u_{i} > 0 ,\\
        \end{array}
 \right. \text{in } \mathbb{R}^n,
\end{align}
where $i=1,\ldots,L$.
The corresponding Dirichlet problem is
\begin{align}\label{systemPotentialDirichlet}
  \left\{\begin{array}{cl}
      	-\Delta u_{i} = \frac{\partial F}{\partial u_i} & \text{in } B, \\
	       u_{i} > 0             & \text{in } B, \\
	       u_i = 0   &\text{on } \partial B,
        \end{array}
 \right.
\end{align}
where $i=1,\ldots,L$, and $B=B_R(0)\subset \mathbb{R}^n$ for any $R>0$.

In \cite{PS-Ind}, Pucci and Serrin have showed that for a general variational problem,
\begin{align*}
    \int_{\Omega}\mathcal{F}(x,u,Du)dx = 0,
\end{align*}
there exists Poho\v{z}aev type of identity that can give
a sufficient condition on the nonexistence of solution to
Dirichlet problem on a bounded star-shaped domain.
In this case, the system \eqref{systemPotential} corresponds to
a vector-valued extremal of the variational problem with
\begin{align*}
    \mathcal{F}(x,u,Du) = \frac{1}{2} \sum_{k=1}^{L} |p^k|^2 - F(u),
\end{align*}
where $p^k=Du_k$. So, theorem 6 in \cite{PS-Ind} leads to the following result (we will also give a simple proof for completion),
\begin{lemma}\label{PucciSerrin}
If for $u\neq 0$,
\begin{align}\label{potentialCondition}
    \frac{n-2}{2}u^k \frac{\partial F}{\partial u^k} - nF(u) > 0,
\end{align}
then system \eqref{systemPotentialDirichlet} with $F(0)=0$ admits no nontrivial solution $u$.
\end{lemma}
Proof. Suppose there exists a solution $u=(u_1,\cdots, u_L)$. 
By a calculation the same as the first step of the proof of Lemma \ref{signChangingDirichletNonEx} we get identities similar to \eqref{pohozaevSingle}, for $ i=1,\cdots,L$,
\begin{align*}
	\int_B \Delta u_i (x\cdot \nabla u_i)  dx &- \int_B  \frac{n-2}{2}|\nabla u_i|^2 dx = \frac 1 2 \int_{\partial B} x\cdot \nu |\nabla u_i|^2 d\sigma>0. 
\end{align*}
Sum the LHS of the above identities up and get
\begin{align*}
 0 &< \int_B -\frac{\partial F}{\partial u_i} (x\cdot \nabla u_i) + \frac{n-2}{2} u_i \Delta u_i  dx  \\
 	&=  \int_B -x\cdot \nabla F(u(x)) - \frac{n-2}{2} u_i \frac{\partial F}{\partial u_i} dx\\
 	&= \int_B n F- \frac{n-2}{2} u_i \frac{\partial F}{\partial u_i} dx,
\end{align*}
which contradicts to \eqref{potentialCondition}.
$\Box$

Therefore, it follows from Theorem \ref{existencetheorem} and Lemma \ref{PucciSerrin} that
\begin{theorem}\label{conservativeEx}
For system \eqref{systemPotential}, if $f$ satisfies the assumptions \eqref{decayAssumption} and \eqref{ControlInequality}, and additionally if $f=\nabla F$, where $F(0)=0$, and $F$ satisfies \eqref{potentialCondition} for $u\neq 0$, then system \eqref{systemPotential} admits a radially symmetric solution of class $C^{2}(\mathbb{R}^{n})$.
\end{theorem}

Here is an example system of {\bf{Type I}}:
Let the potential function be
\begin{align}
	F(u,v)=-(u-v)^2+v^{p-1}u+u^{p-1}v, \text{with } p\geq \frac{2n}{n-2},
\end{align}
 and then
\begin{align}\label{example3}
  \left\{\begin{array}{ll}
    -\Delta u = F_u = -2(u-v) + v^{p-1} + (p-1)u^{p-2}v,  \\
    -\Delta v = F_v = 2(u-v) + (p-1)v^{p-2}u + u^{p-1}, \\
	\ \ u,v > 0, \\
        \end{array}
\right. \text{in } \mathbb{R}^n.
\end{align}
We can verify that $F$ satisfies the condition of Theorem \ref{conservativeEx}. $F(0,0)=0$, and $F_u+F_v \geq 0$ (\eqref{decayAssumption} is satisfied). Let $u=0$ then $F_u = 2v+v^{p-1}\geq |F_v|=2v$, similarly if $v=0$, then $F_v =2u+u^{p-1}\geq |F_u|=2u$, so we can see that \eqref{ControlInequality} is satisfied in a neighborhood. Last, direct computation shows that
\begin{align*}
	\frac{n-2}{2}(uF_u+vF_v)-nF &= \frac{n-2}{2}(-2(u-v)^2 + puv^{p-1} + pu^{p-1}v) -n(-(u-v)^2+v^{p-1}u+u^{p-1}v) \\
	&\geq 2(u-v)^2 >0, 
\end{align*}
so \eqref{potentialCondition} is satisfied. Also, notice that $F_u$ and $F_v$ are sign-changing functions, for example, let $v=0$, then $F_u<0$ and let $u=0$ then $F_u>0$.

{\bf{Type II.}}
In the following example, we give another class of systems with potential type of source term. 

Consider the following system, in $\mathbb{R}^n$,
\begin{align}\label{systemPotential2}
  \left\{\begin{array}{cl}
      	-\Delta u = f_1=\frac{\partial F}{\partial v},  \\
      	-\Delta v = f_2 = \frac{\partial F}{\partial u}, \\
	       u,v > 0 ,\\
        \end{array}
 \right. 
\end{align}
and corresponding Dirichlet problem
\begin{align}\label{systemPotentialDirichlet2}
  \left\{\begin{array}{cl}
      	-\Delta u = \frac{\partial F}{\partial v}, \text{in }  B \\
      	-\Delta v = \frac{\partial F}{\partial u}, \text{in }  B \\
	       u,v = 0 , \text{on }  \partial B\\
        \end{array}
 \right. 
\end{align}
where $B=B_R(0)\subset \mathbb{R}^n$ for any $R>0$.
Similarly we have
\begin{lemma}\label{PucciSerrin2}
If for $(u,v)\neq (0,0)$,
\begin{align}\label{potentialCondition2}
    \frac{n-2}{2}\left( u \frac{\partial F}{\partial u} + v \frac{\partial F}{\partial v}\right)  - nF(u,v) > 0,
\end{align}
then system \eqref{systemPotentialDirichlet2} with $F(0)=0$ admits no nontrivial solution $u$.
\end{lemma}
Then it follows from Theorem \ref{existencetheorem} and Lemma \ref{PucciSerrin2} that
\begin{theorem}\label{conservativeEx2}
For system \eqref{systemPotential2}, if $f_1=F_v$ and $f_2=F_u$ satisfies the assumptions \eqref{decayAssumption} and \eqref{ControlInequality}, and additionally if $F(0)=0$, and $F$ satisfies \eqref{potentialCondition2} for $(u,v)\neq (0,0)$, then system \eqref{systemPotential2} admits a radially symmetric solution of class $C^{2}(\mathbb{R}^{n})$.
\end{theorem}

{\bf{Proof of Lemma \ref{PucciSerrin2}.}}
Suppose there exists a solution $(u,v)$. 
Then we follow the proof of Lemma \ref{signChangingDirichletNonEx}. The first step is the same, and we merge the source terms to \eqref{pohozaevCross} and get,
\begin{align*}
0 &<\int_B \Delta u(x\cdot \nabla v) + \Delta v(x\cdot\nabla u) - (n-2)\nabla u\cdot \nabla v dx \\
	&= \int_B -F_v(x\cdot \nabla v) - F_u(x\cdot\nabla u)+ (n-2)(\theta v\Delta u + (1-\theta) u\Delta v)dx  \\
	&= \int_B -x\cdot \nabla F -(n-2)(\theta vF_v + (1-\theta) uF_u) dx\\
	&= \int_B nF-(n-2)(\theta vF_v + (1-\theta) uF_u) dx,
\end{align*}
which contradicts to \eqref{potentialCondition2} when taking $\theta=\frac 1 2$.
$\Box$

Here is an example system of {\bf{Type II}}:
Let the potential function be
\begin{align}
F(u,v)=-(u-v)^2+u^p+v^p, \text{with } p\geq \frac{2n}{n-2},
\end{align}
and
\begin{align}
	  \left\{\begin{array}{ll}
	    -\Delta u = F_v = -2(u-v) + pv^{p-1},  \\
	    -\Delta v = F_u = 2(u-v) +  pu^{p-1}, \\
		\ \ u,v > 0, \\
	        \end{array}
	\right. \text{in } \mathbb{R}^n.
\end{align}
We verify that $F$ satisfies conditions in Theorem \ref{conservativeEx2}. $F(0)=0$, and $F_u+F_v \geq 0$ (\eqref{decayAssumption} is satisfied). Let $u=0$ then $F_v = 2v+v^{p-1}\geq |F_u|=2v$, similarly if $v=0$, then $F_u =2u+u^{p-1}\geq |F_v|=2u$, so we can see that \eqref{ControlInequality} is satisfied in a neighborhood. Last, direct computation shows that
\begin{align*}
	\frac{n-2}{2}(uF_u+vF_v)-nF &= \frac{n-2}{2}(-2(u-v)^2 + pu^p + pv^p) -n(-(u-v)^2+v^p+u^p) \\
	&\geq 2(u-v)^2 >0, 
\end{align*}
so \eqref{potentialCondition2} is satisfied. Also, notice that $F_u$ and $F_v$ are sign-changing functions, for example, let $v=0$, then there exists sufficiently small $u$ such that $F_u<0$, and let $u=0$ then there exists sufficiently small $v$ such that $F_v>0$.

\end{document}